\documentclass[11pt,a4paper]{amsart}
\usepackage[utf8]{inputenc}
\usepackage{amsmath}
\usepackage{amsfonts}
\usepackage{amssymb, amscd, amsthm}
\usepackage{graphicx}
\usepackage{fullpage}
\usepackage[all]{xy}
\usepackage{enumitem, hyperref, comment}

\newtheorem{theorem}{Theorem}

\newtheorem{corollary}[theorem]{Corollary}
\newtheorem{conjecture}[theorem]{Conjecture}
\newtheorem{proposition}[theorem]{Proposition}

\title{An introduction to the Kepler-Heisenberg problem}
\author{Corey Shanbrom}
\address{California State University, Sacramento, 6000 J St., Sacramento, CA 95819, USA}
\email{corey.shanbrom@csus.edu}
\date{\today}

\begin{document}
\maketitle

\section{Introduction}

Here we provide an overview of what is known, and what is not known, about an interesting dynamical system known as the Kepler-Heisenberg problem.  The main idea is to pose a version of the classical Kepler problem of planetary motion, but in a sub-Riemannian setting.  The result is system which is surprisingly rich and beautiful, mysterious in some ways but tame in others, offering a substantial number of questions which seem non-trivial yet tractable.

\subsection{History and motivation}
The curved Kepler problem, typically set in a Riemannian manifold of constant curvature, is an active area with a fascinating history.  In particular, while the co-discovery of hyperbolic space by Lobachevsky and Bolyai is famous mathematical lore, both men also independently attempted to pose the Kepler problem in this new geometry around 1835.  Dirichlet and Schering worked on this problem shortly after.  In 1873, Lipschitz posed and partially solved the Kepler problem on the 3-sphere, as did Killing (with a different notion of gravity) in 1885.  In 1902, Liebmann proved that orbits are conics in both hyperbolic and and spherical 3-space.  
There is a parallel history of studying the Kepler problem on surfaces as well, and today there is active work on $n$-body problems in more general curved spaces.
Details can be found in \cite{Diacu} and the references therein (while this paper is published, the arXiv version contains a much more elaborate history).

The original motivation for studying the Kepler-Heisenberg problem was to extend this idea to \textit{sub-Riemannian} manifolds, which generalize Riemannian ones.  In particular, the three dimensional Heisenberg group is the simplest non-trivial sub-Riemannian space and the natural place to begin such an investigation.  
While no knowledge of sub-Riemannian geometry is needed for this article, the interested reader can see \cite{Tour} for an introduction to this subject in general, and details of the Heisenberg group in particular.  
As seen below, the resulting Kepler-Heisenberg system retains certain aspects of the classical Kepler problem in Euclidean space, while reflecting the unique geometry of the Heisenberg group.  We are thus left with obvious questions about generalizing both the setting to other sub-Riemannian spaces and the force to other classical mechanical systems.  See Section \ref{open}.

We can provide an alternative motivation as well.  The classical Kepler problem enjoys the famous three Kepler's Laws of Planetary Motion.  Do these hold in curved spaces?  For spherical and hyperbolic 3-space, the first two laws hold but the third fails.  Why?  Kepler's third law is a consequence of homogeneity of the gravitational potential with respect to Euclidean dilations.  But spherical and hyperbolic spaces do not admit such dilations.  Moreover, by a result of Gromov (\cite{Gromov}), there are no homogeneous Riemannian manifolds which admit dilations besides Euclidean spaces.  Thus, if one wants to pose a curved Kepler problem in which Kepler's third law holds, one must leave the Riemannian realm.  Certain important sub-Riemannian spaces known as Carnot groups do admit dilations.  The Heisenberg group is the simplest, and we will see below that a version of Kepler's third law does indeed hold there.

\subsection{The problem}
We treat the Heisenberg group $\mathcal H$ as the smooth manifold $\mathbb R^3$ equipped with a certain non-Euclidean metric.  
Take the usual global coordinates $(x, y, z)$ and define two vector fields
$$X=\frac{\partial}{\partial x}-\frac{1}{2}y\frac{\partial}{\partial z}\quad \text{and}\  \quad Y=\frac{\partial}{\partial y}+\frac{1}{2}x\frac{\partial}{\partial z}.
$$
These span a two-dimensional bracket-generating distribution $\mathcal D$.  The Kepler-Heisenberg problem is roughly the classical Kepler problem in $\mathbb R^3$ with the constraint that all orbits are tangent to $\mathcal D$; such curves are called \textit{horizontal}.  Declaring that $X$ and $Y$ are orthonormal defines an inner product on $\mathcal D$.    Then $\mathcal H$ is now equipped with a sub-Riemannian structure whereby we can measure lengths of horizontal curves and (by the Chow-Rashevskii theorem) distances between points.  The resulting geometry enjoys three interesting properties:
\begin{itemize}
\item the length of a horizontal curve is equal to the Euclidean length of its projection to the $xy$-plane ($ds^2=(dx^2 +dy^2)|_{\mathcal D}$),
\item the $z$-coordinate of a horizontal curve grows like the area traced out by its projection (Stokes' theorem),
\item geodesics are helices projecting to circles or lines (dual isoperimetric problem).
\end{itemize}

Our configuration space is $\mathcal H$ equipped with this structure.  Here we pose the Kepler-Heisenberg problem as a Hamiltonian system, taking $T^*\mathcal H$ with coordinates $(x, y, z, p_x, p_y, p_z)$ as phase space.  Consider the dual momenta to the vector fields $X$ and $Y$ given by
$$P_X=p_x-\tfrac{1}{2}yp_z \quad \text{and}\  \quad P_Y=p_y+\tfrac{1}{2}xp_z.
$$
Then, as in the Riemannian case, we take our kinetic energy to be 
$K=\tfrac{1}{2}(P_X^2+P_Y^2)$.  This is induced by our metric and generates geodesic flow.  
Choosing the gravitational potential is more subtle.  In the classical case, the potential $1/r$ has many important properties.  It is most obviously the reciprocal of distance from the origin, but for most (sub-)Riemannian manifolds there is no explicit expression for such distance, so this characterization is inadequate.  Instead, we recognize $1/r$ as (a multiple of) the fundamental solution to the Laplacian on $\mathbb R^3$, and choose our potential accordingly.  Thus we choose the Kepler-Heisenberg potential as $U=-\frac{1}{8\pi}\Big((x^2+y^2)^2+16z^2\Big)^{-1/2}$, which is the fundamental solution to the Heisenberg sub-Laplacian $\Delta_{\mathcal{H}}=X^2 +Y^2$ according to \cite{Folland}.  Note the singularity at the origin, which represents our planet colliding with the sun.

Thus we obtain our Hamiltonian 
$$H=\underbrace{\tfrac{1}{2}((p_x-\tfrac{1}{2}yp_z)^2+(p_y+\tfrac{1}{2}xp_z)^2)}_K\underbrace{-\frac{1}{8\pi\sqrt{(x^2+y^2)^2+{16}z^2}}}_{U}.$$
The \textit{Kepler-Heisenberg problem} is to analyze and eventually solve Hamilton's equations for $H$.  	

\section{Results}
We first give some basic properties of the system, most of which can be found in \cite{MS}.  An important tool will be the Heisenberg dilation
$(x, y, z)\mapsto  (\lambda x, \lambda y, \lambda^2 z)$
for $\lambda>0$.

\begin{proposition}
The Kepler-Heisenberg system has two integrals of motion: the total energy $H$ and the angular momentum  $p_{\theta}=xp_y-yp_x$ .  The dilational momentum $J=xp_x+yp_y+2zp_z$ satisfies $\dot J = 2H$.  The dynamics are therefore integrable on the codimension one submanifold $\{H=0\}$, where $J$ is also conserved.  Periodic orbits must satisfy $H=J=0$.
\end{proposition}
These are the known symmetries -- there may be others.  So we do not know whether this system is integrable or not.  But we consider it `at least mostly integrable'.

Ostensibly, we should be able to integrate the zero-energy dynamics.  However, the dilational action is non-compact, so the standard machinery of action-angle variables does not immediately apply.  (While there is a generalized theory, as in \cite{Fiorani}, it is non-constructive.)  Our invariant submanifolds can be visualized as cones of 2-tori,  $\mathbb T^2 \times \mathbb R^+$, with the cone point at $0 \in \mathbb R^+$ representing the singularity at sun.  If an orbit's constant dilational momentum is positive (resp. negative) then it winds away (resp. towards) the cone point; if $J=0$ then the orbit stays on a compact $\mathbb T^2$ and the motion is quasi-periodic or periodic.  In \cite{MS},  the integration of the $H=0$ system is reduced to parametrizing a family of degree 6 algebraic curves, but there is no general method for obtaining such a parametrization.

As mentioned in the Introduction, we also have a version of Kepler's third law:
$$T^2 = ka^4$$
where $T$ represents a orbit's period, $a$ is its size, and $k$ is a universal constant.  This is a consequence of our Hamiltonian being homogeneous of degree $-2$ with respect to the Heisenberg dilation.  Note that we do have Kepler's second law as well, since it is equivalent to the conservation of angular momentum.  Kepler's first law, however, fails in the Heisenberg group (recall that it held in spherical and hyperbolic three-space, where the third law failed).  We do have periodic orbits (see below), but they are not conic sections.

Finally, we mention a few remaining small but interesting results.  All negative energy solutions are bounded.  The only solutions contained in the $xy$-plane are lines through the origin.  On the invariant submanifold $\{z=p_z=p_{\theta}=0\}$ we have a classical central force problem in the plane, and the radial distance $r(t)$ traces out a hyperbola, ellipse, or parabola if the energy is positive, negative, or zero, respectively.  Lastly, there is one very strange family of orbits: each sits at a constant position on the $z$-axis while its momenta $p_z$ grows linearly in time.

\subsection{Periodic orbits} 
\begin{figure}
\centering
\begin{tabular}{cc}
\includegraphics[width=.25\textwidth]{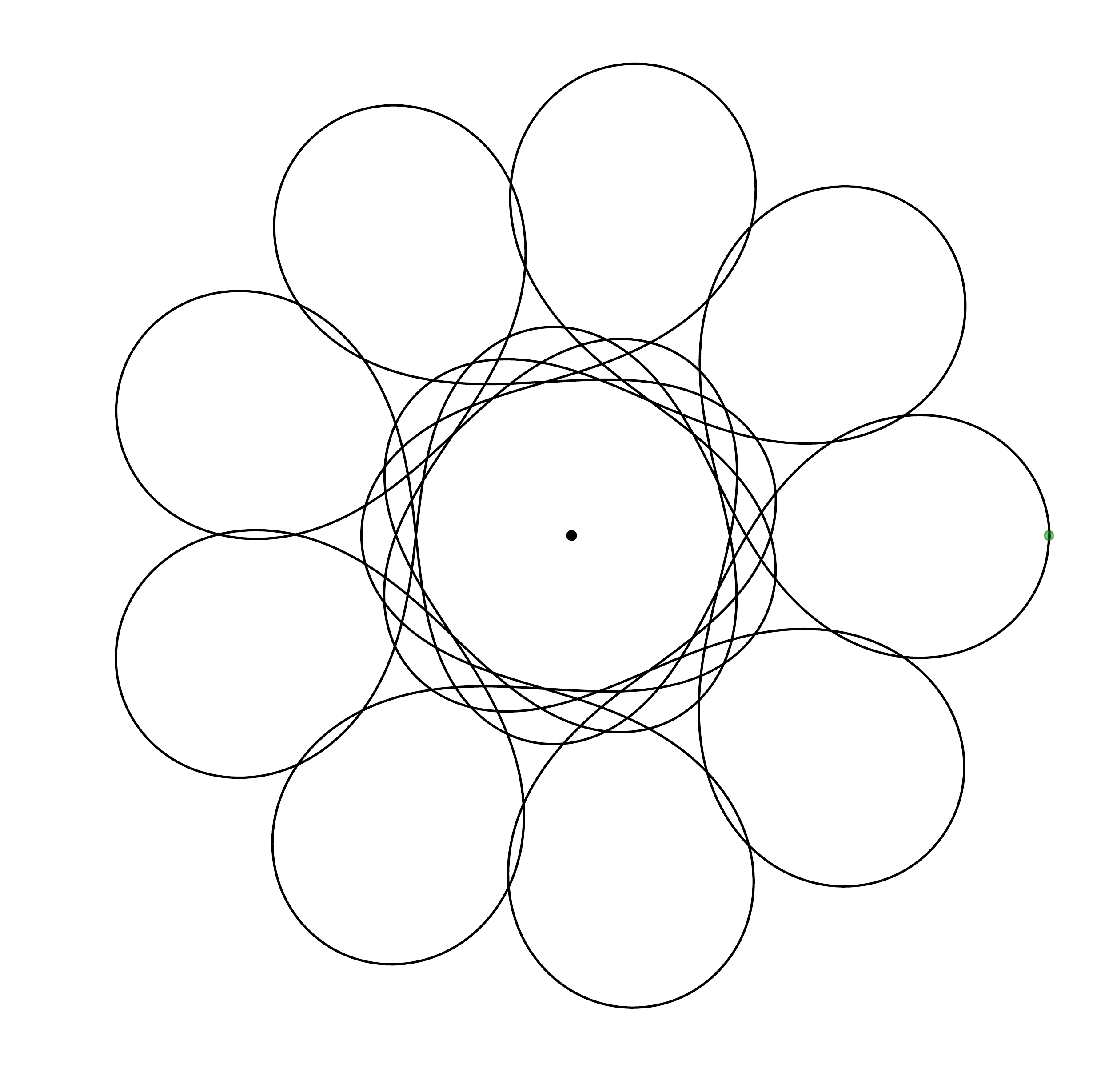}
\includegraphics[width=.25\textwidth]{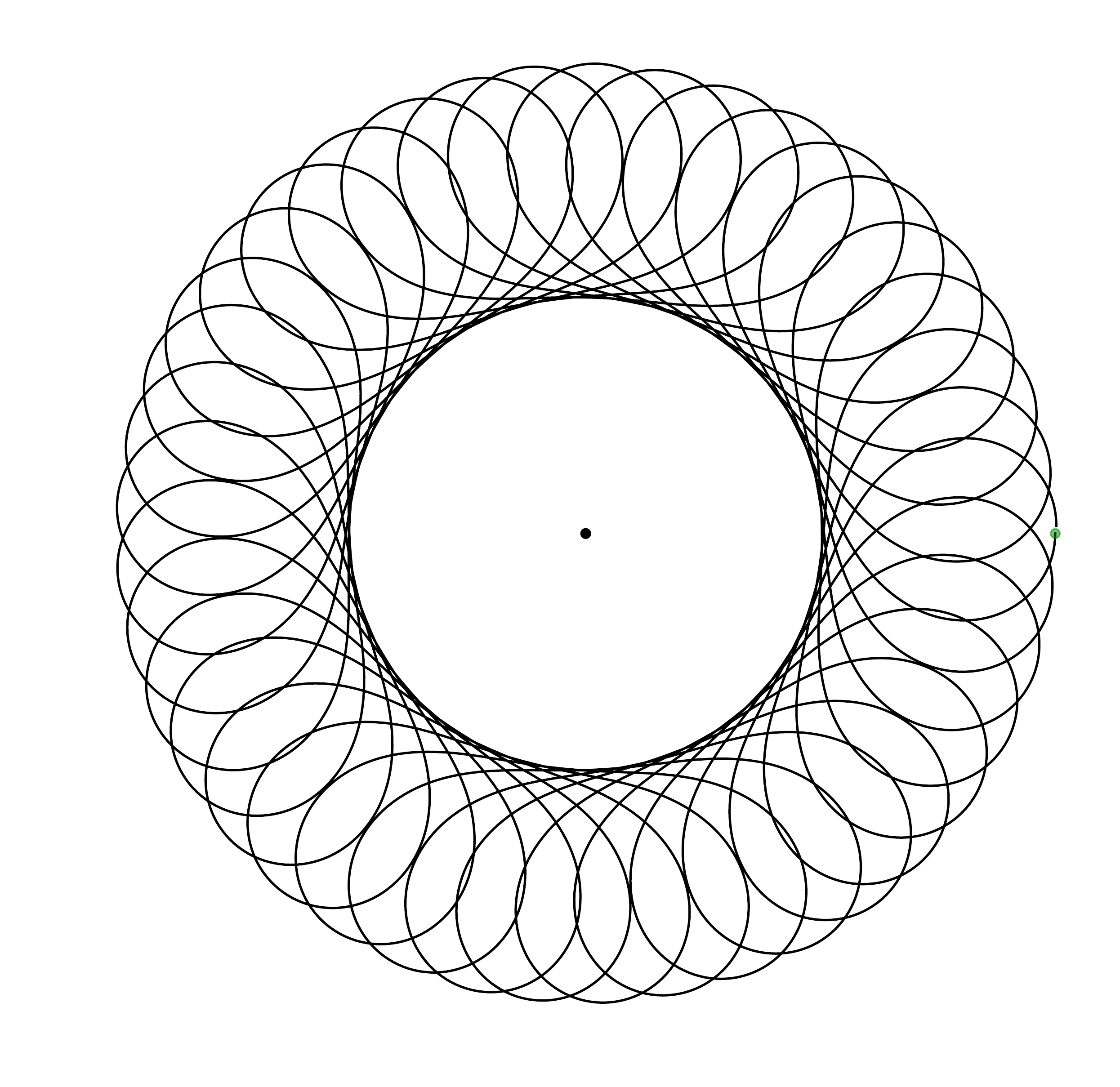}
\\
\includegraphics[width=.25\textwidth]{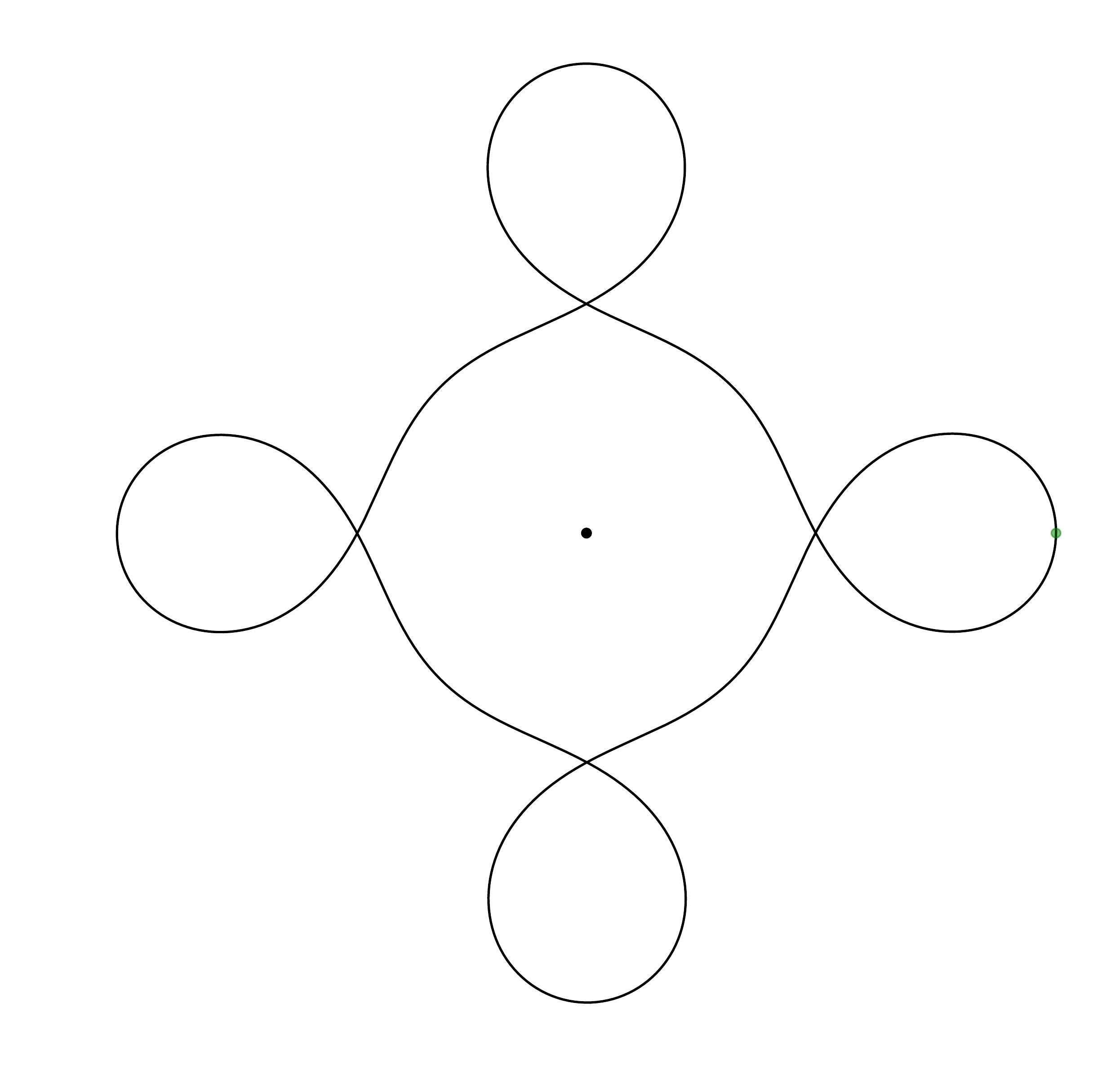}
\includegraphics[width=.25\textwidth]{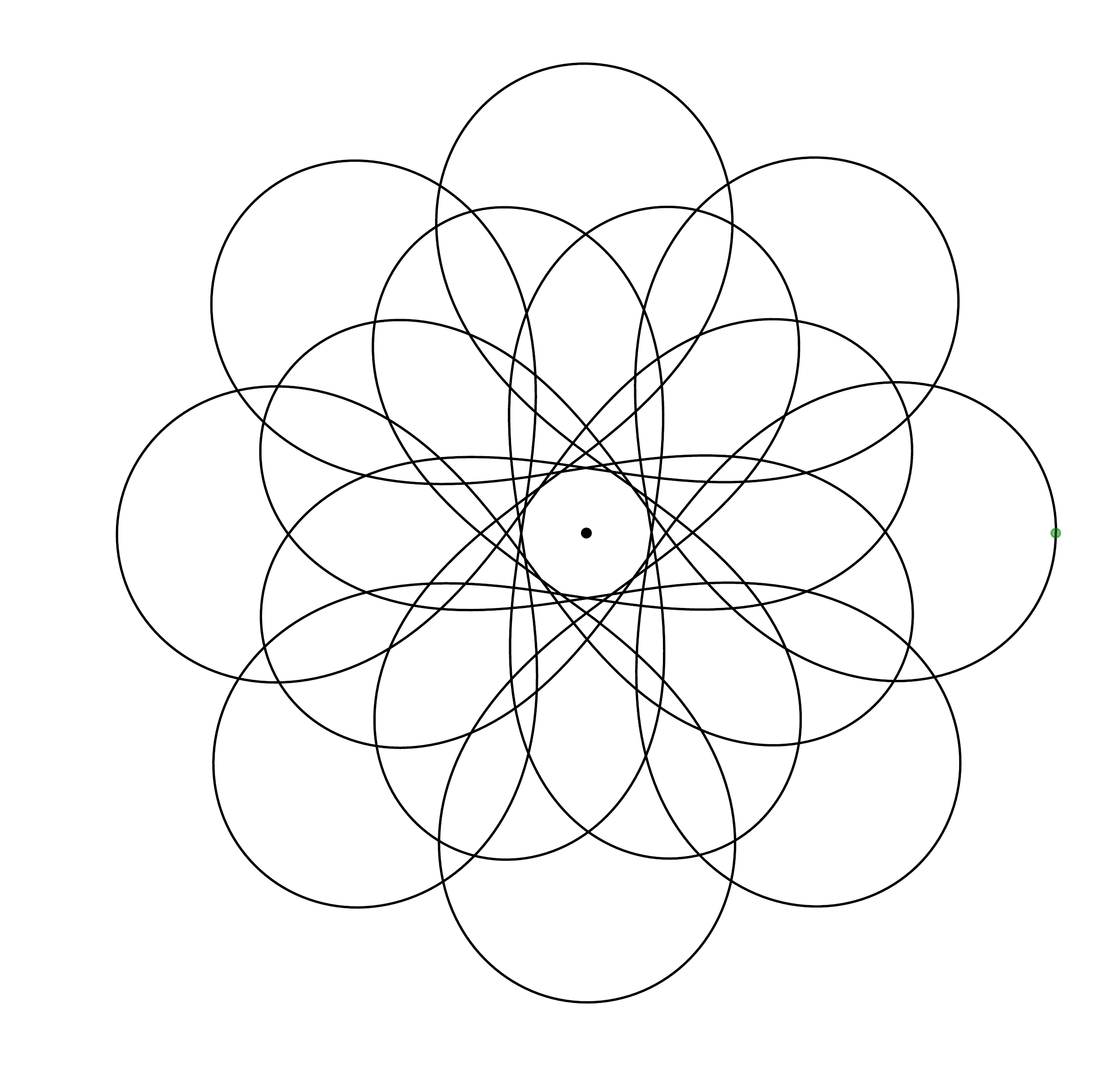}
\end{tabular}
\caption{Periodic orbits with rotation numbers 5/9, 6/41, 1/4, and 7/8.  In all figures we show the  projection of the orbit to the $xy$-plane.; the $z$-coordinate grows like the area traced out.
} \label{periodic fig 1}
\end{figure}

\begin{figure}
\begin{centering}
\includegraphics[width=0.23\textwidth]{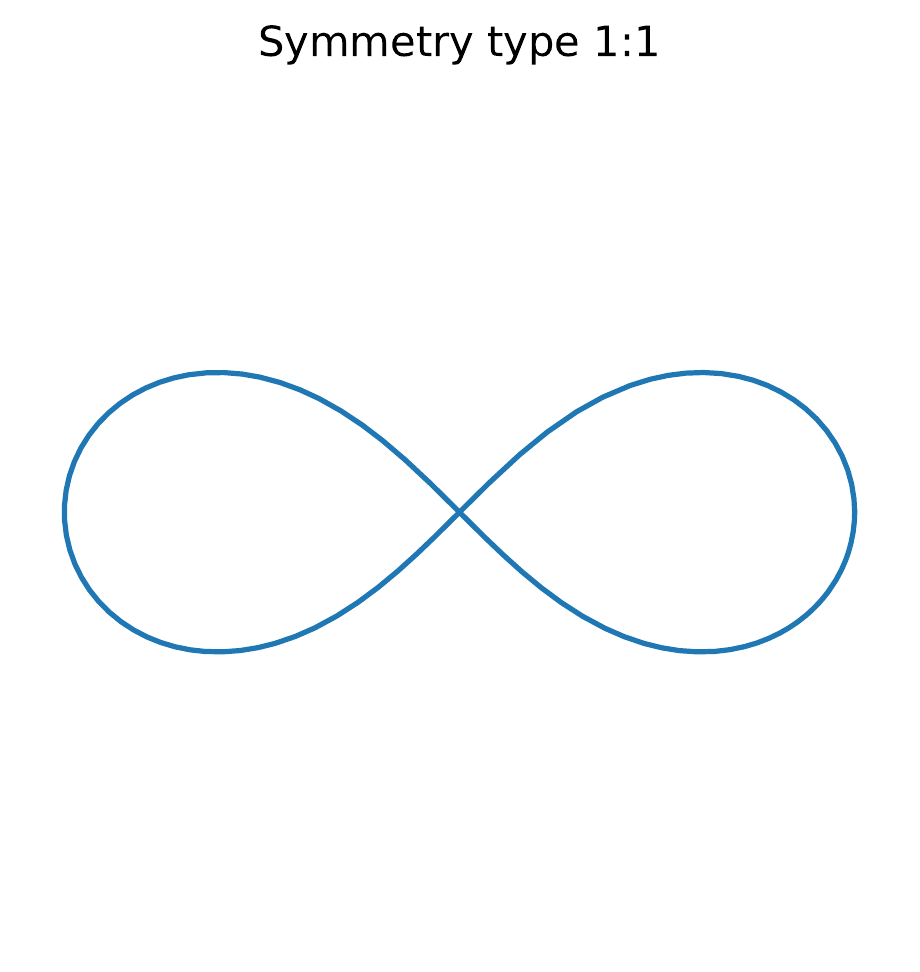}
\includegraphics[width=0.23\textwidth]{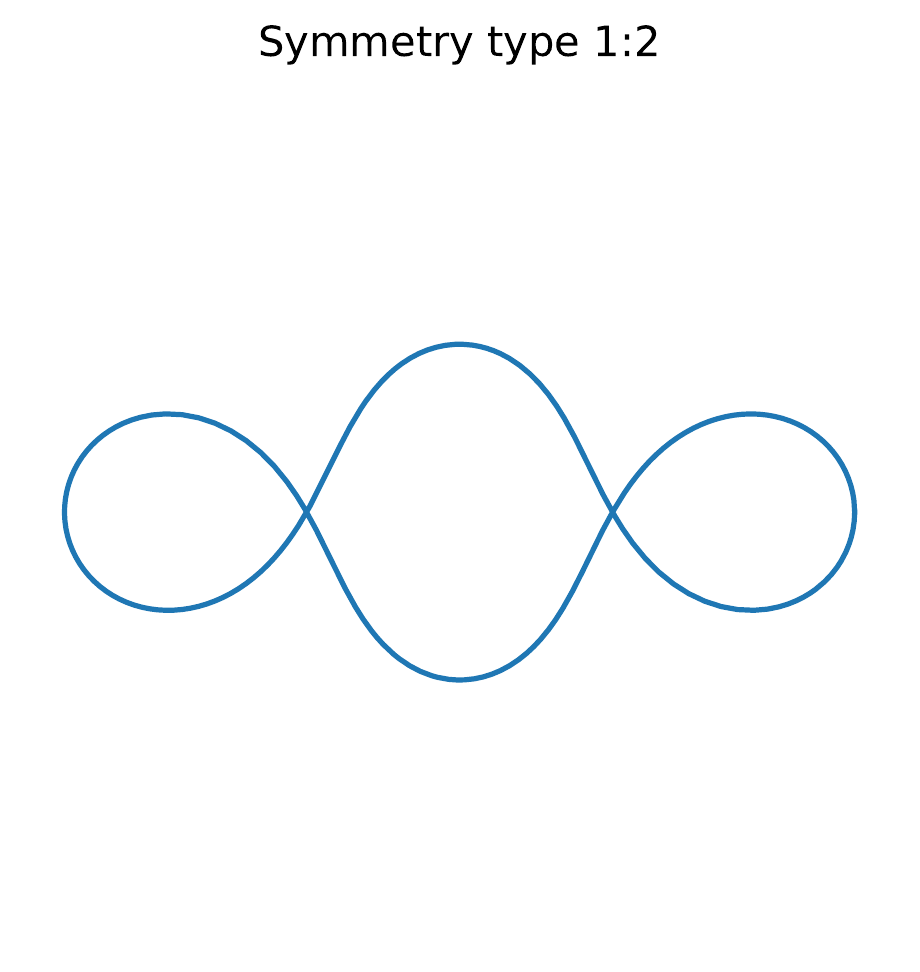}
\par\end{centering}
\begin{centering}
\includegraphics[width=0.23\textwidth]{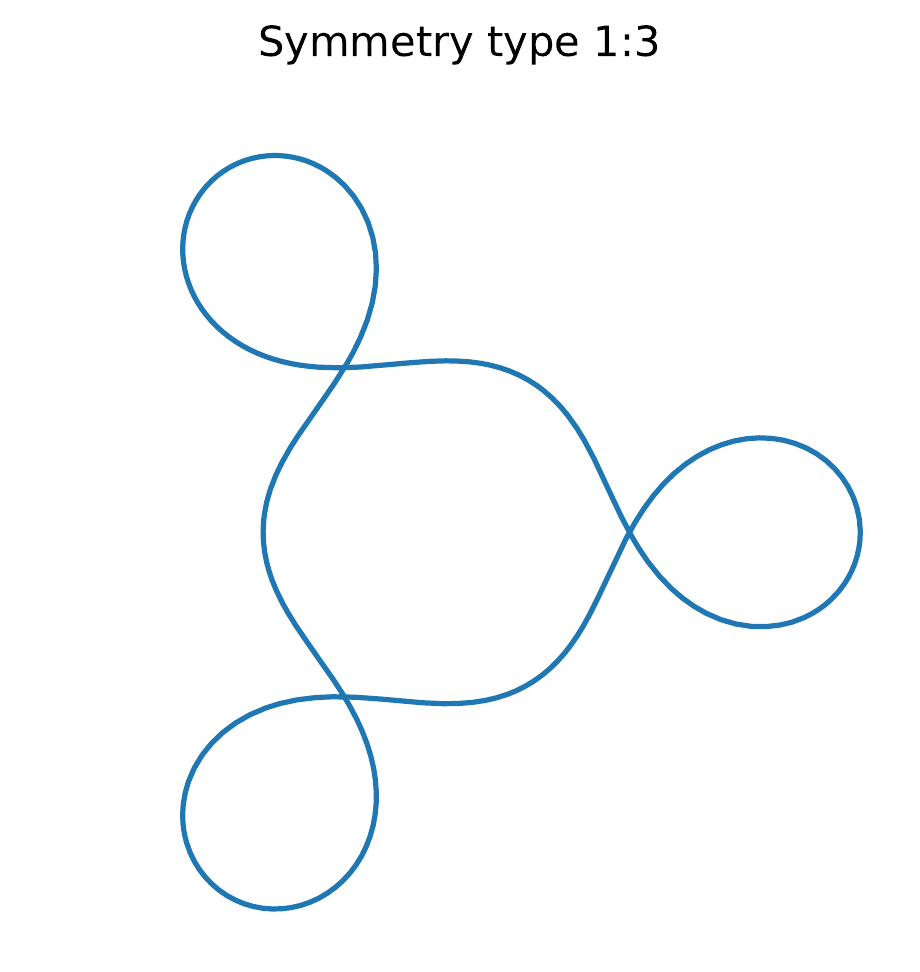}
\includegraphics[width=0.23\textwidth]{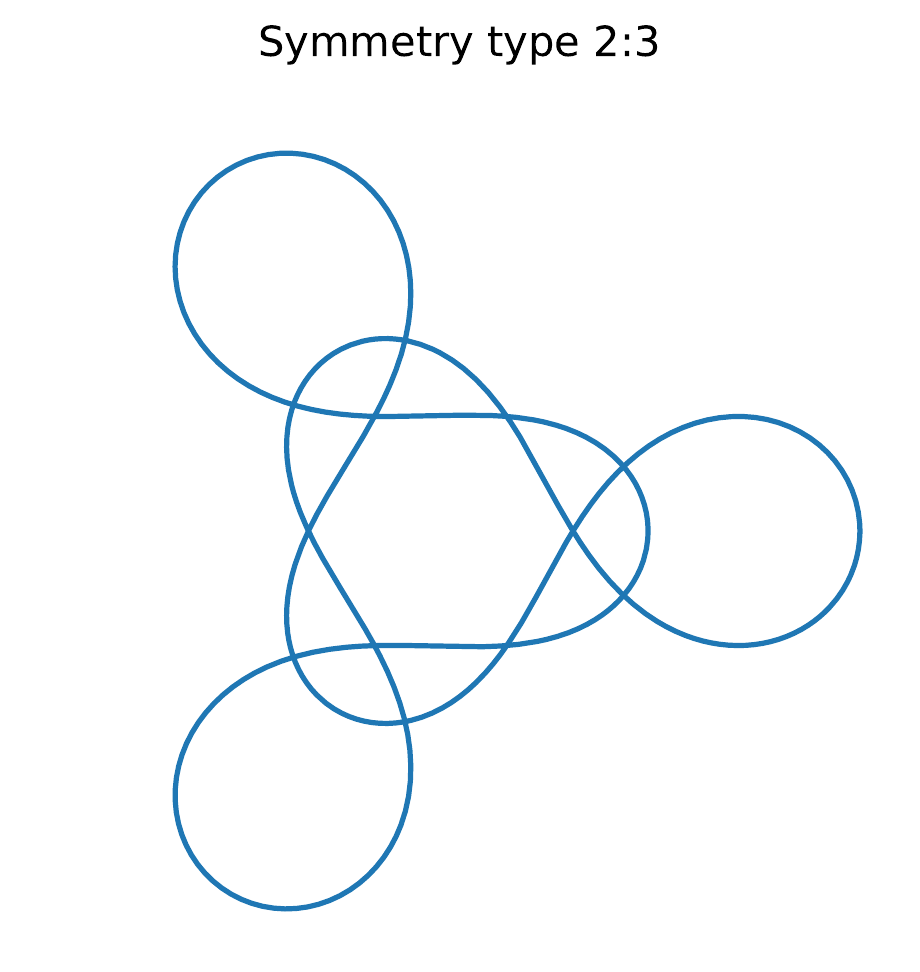}
\includegraphics[width=0.23\textwidth]{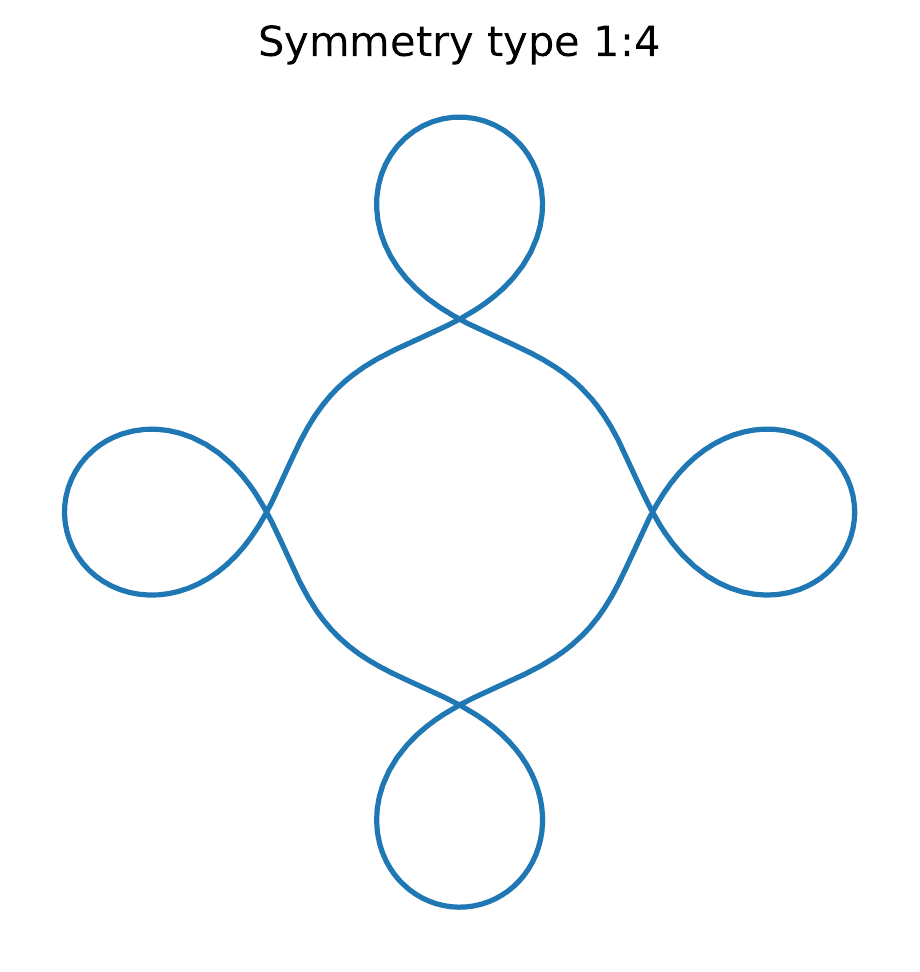}
\includegraphics[width=0.23\textwidth]{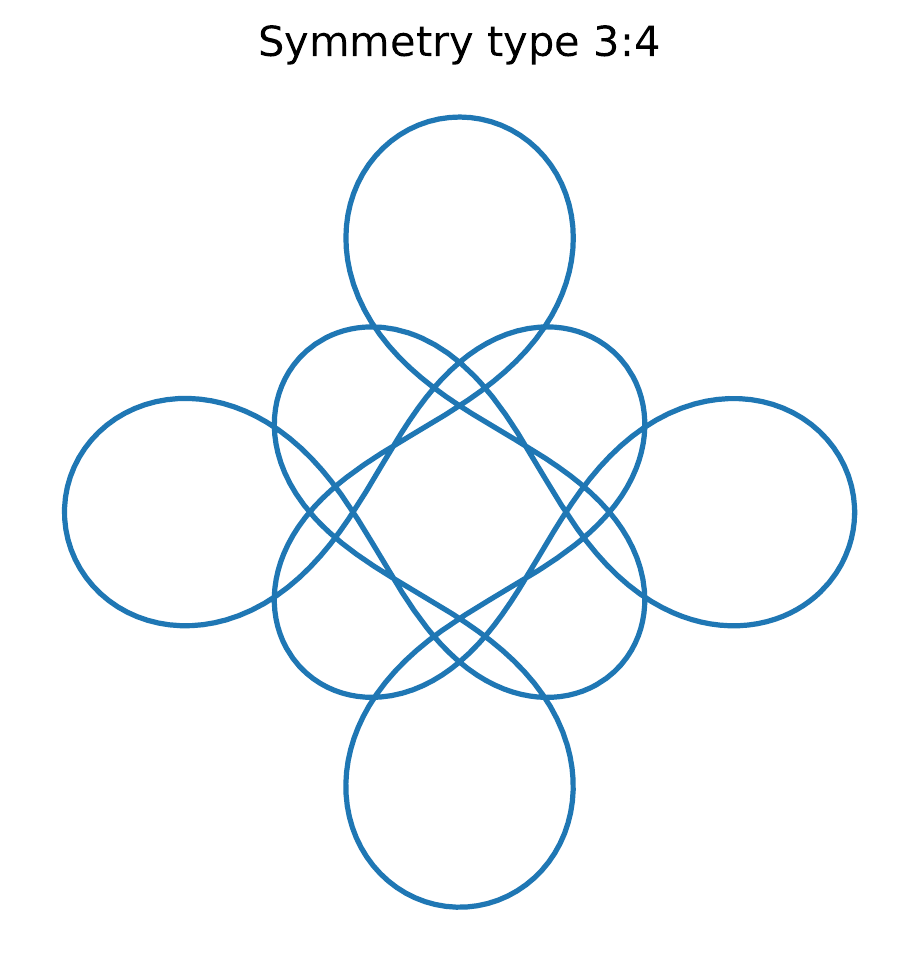}
\par\end{centering}
\begin{centering}
\includegraphics[width=0.23\textwidth]{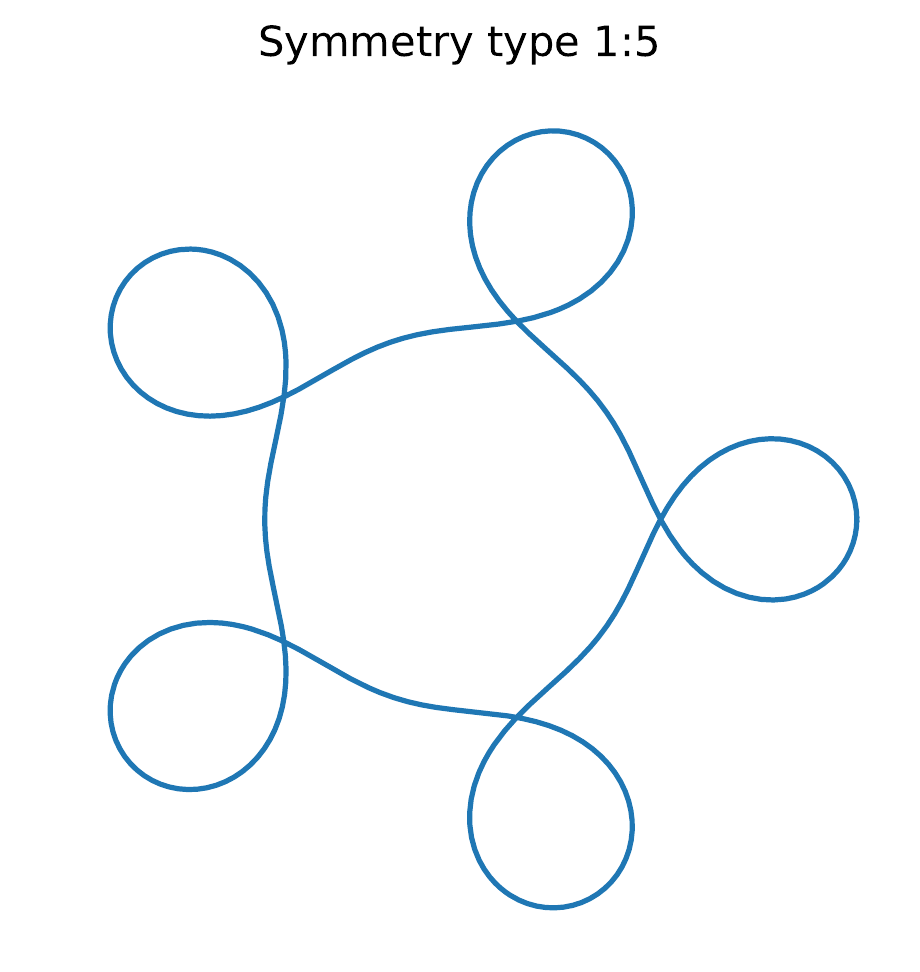}
\includegraphics[width=0.23\textwidth]{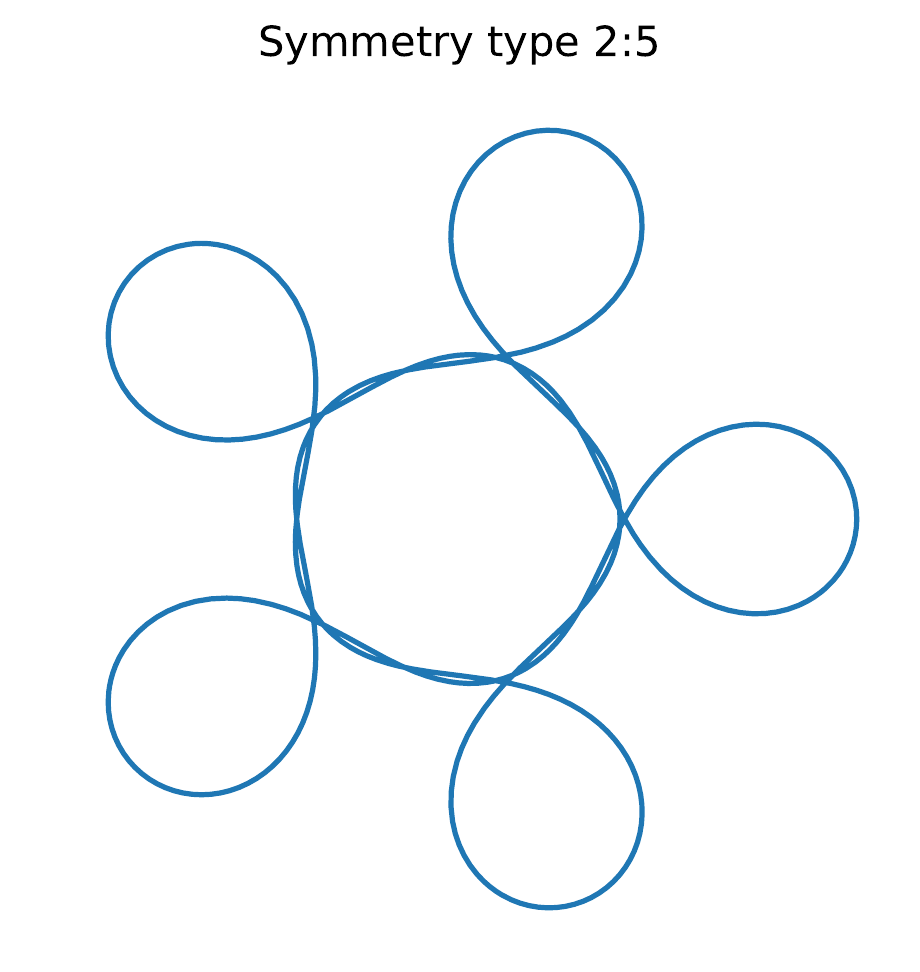}
\includegraphics[width=0.23\textwidth]{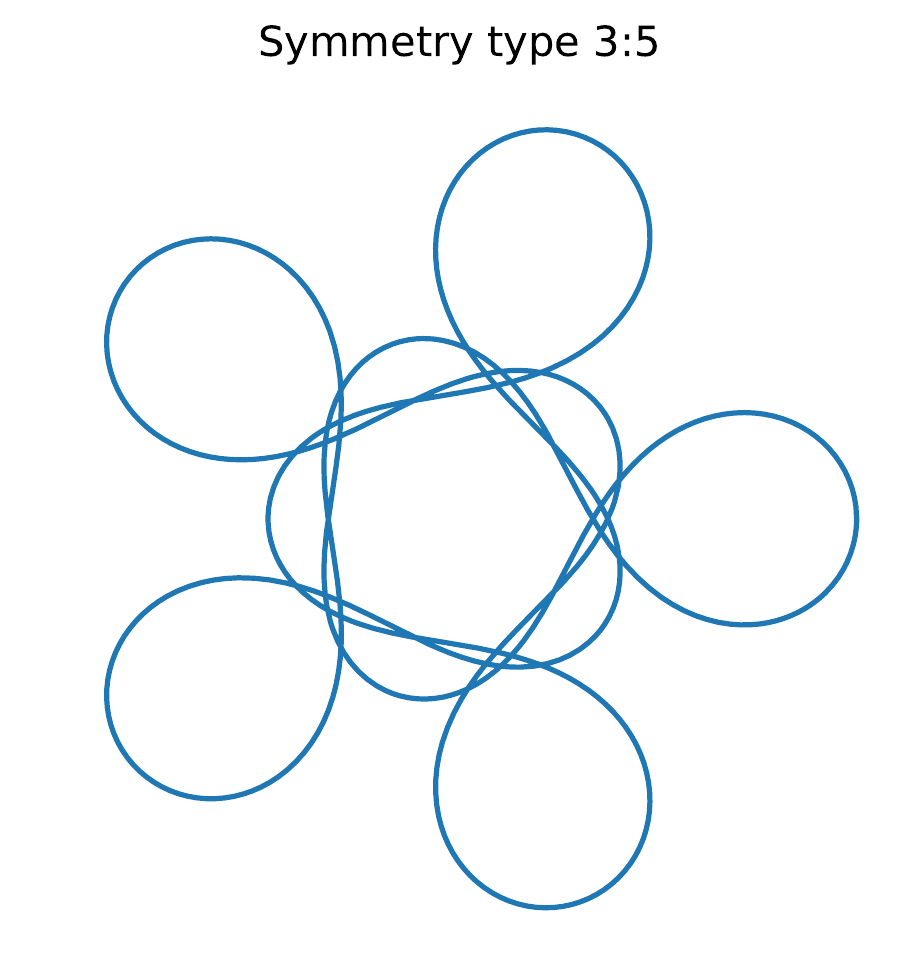}
\includegraphics[width=0.23\textwidth]{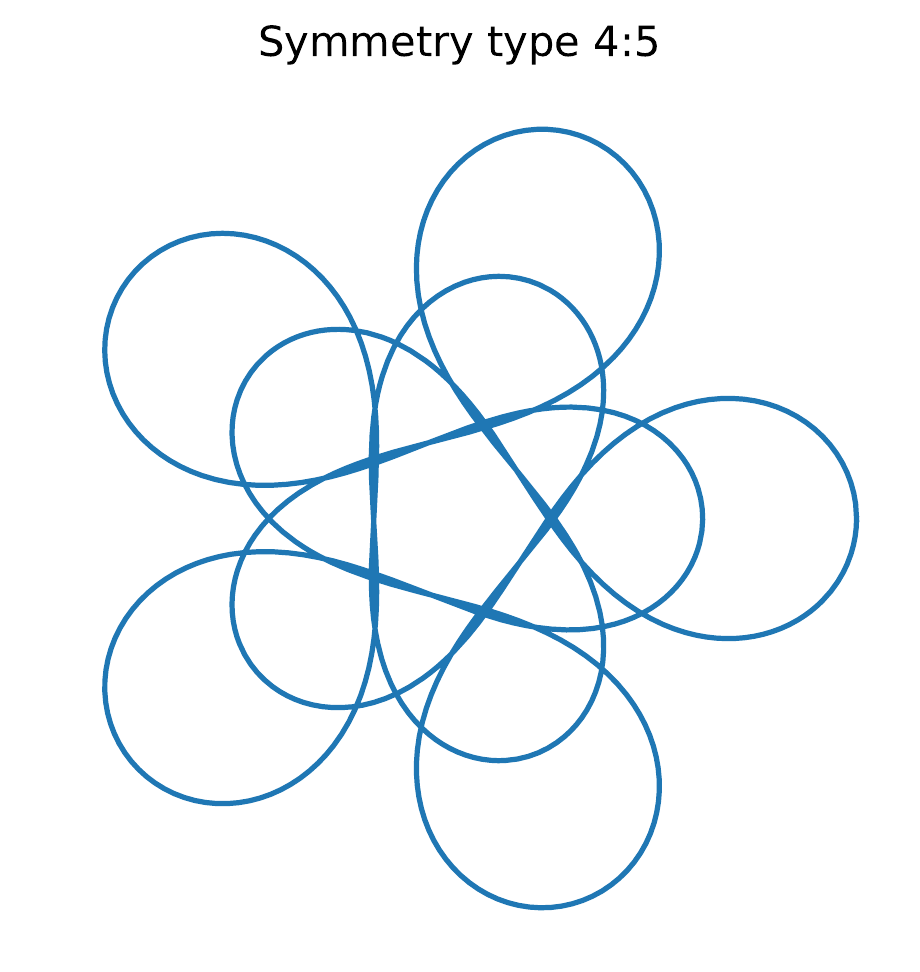}
\par\end{centering}
\caption{All periodic orbits with $k\leq 5$.}\label{gallery}
\end{figure}

The existence of periodic orbits was established in \cite{CS}.  However, the proof followed the direct method in the calculus of variations and provided very little information about such orbits beyond their existence.  Numerical methods in \cite{DS} allowed the discovery of a very rich and beautiful symmetry structure.
\begin{theorem}
For any rational $j/k \in (0,1]$, there is a periodic Kepler-Heisenberg orbit with rational rotation number $j/k$.  
\end{theorem}
Rational rotation numbers are treated in \cite{Katok, Tabachnikov}; here it simply means that
after one $k$th of the period, the planet has rotated $2\pi \cdot \frac{j}{k}$ radians counterclockwise about the $z$-axis.  See Figures \ref{periodic fig 1} and \ref{gallery}.  All such orbits possess $k$-fold rotational symmetry about the $z$-axis, except the orbit with rotation number 1/1 (its projection does, but not $z(t)$).  

This symmetry structure was uncovered via numerical explorations. We strongly
 encourage other researchers to use our codebase \cite{Github}.  All code is hosted publicly as free, open-source software; explicit instructions to reproduce each plot are provided there, with all relevant documentation and licensing.  All results in \cite{DS} are thus 100\% available and exactly reproducible.

Our Python program searches for a closed orbit nearby some given initial condition, which may be input directly or randomly generated.  First, the program checks whether this starting initial condition $X_0$ is close enough to a closed orbit.  To do so we use a symplectic integrator and shooting method with objective function \textit{obj} giving the smallest local minimum of distance from $X_0$.
 If $X_0$ is satisfactory, we use a Monte Carlo method to optimize \textit{obj} as follows.  
Randomly pick a point in radially symmetric distribution on an annular neighborhood centered at $X_0$, then symplectically integrate and evaluate \textit{obj}.
 If \textit{obj} does not decrease, try a new point; if \textit{obj} decreases, use this point as new $X_0$.  Repeat this {update} step 1000 times.
 The output is an image, a binary data file, and a human readable summary containing all pertinent information about the orbit (initial conditions, angular momentum, period, etc.).  The program also detects the rotation number using Fourier analysis.
 
Finally, our numerical investigations discovered another somewhat surprising phenomenon that we do not fully understand.   Due to the symmetries and constraints of the system, we can reduce the search for periodic orbits to one dimension, which we can parametrize by angular momentum $p_{\theta}$.  Then as a function of  $p_{\theta}$, the rotation numbers are arranged according to the Farey sequence.  See Table \ref{table}.

\begin{table}[h!]
  \caption{As $p_{\theta}$ grows, the rotation numbers are distributed like the Farey sequence. }
\renewcommand\arraystretch{1.5}
  \begin{tabular}{ | c || c | c | c | c | c | c | c | c | c | c | c | c | }
    \hline
    $\approx p_{\theta}$ & 3.04e--6 & .060 & .071 & .087 & .113 &  .133 & .164 & .199 & .226 & .271 & .307 & .339 \\              
    \hline
    $\frac{j}{k}$ & $\frac{1}{1}$ & $\frac{5}{6}$ & $\frac{4}{5}$ & $\frac{3}{4}$ & $\frac{2}{3}$ & $\frac{3}{5}$ & $\frac{1}{2}$ & $\frac{2}{5}$ & $\frac{1}{3}$ & $\frac{1}{4}$ & $\frac{1}{5}$ & $\frac{1}{6}$ \\ \hline
  \end{tabular}
  \label{table}
\end{table}


\subsection{Self-similarity}
Zero-energy orbits are relatively well behaved and understood.  Almost all exhibit an attractive self-similarity property (\cite{DS2}).  
To the best of our knowledge these solutions represent the first occurrence of this type of self-similarity within 
solutions to ODEs of Hamiltonian type.
\begin{theorem}\label{self sim thm}
Suppose $c
$ is a curve in $T^*\mathcal{H}$ such that:
\begin{itemize}
\item $c$ solves Hamilton's equations
\item $H(c(t))=0$
\item $z$ has consecutive zeros at $t_0, t_1, t_2$
\end{itemize} 
Then $c$ is entirely determined by its value on the fundamental domain $[t_0, t_2]$.   
\end{theorem}

See Figure \ref{self sim fig}, which is explained in detail below. The self-similarity here involves a rotation, a dilation, and a time reparametrization.  The first two assumptions in the theorem are clear, but the third requires some explanation.  Numerical evidence suggests the following.

\begin{conjecture}\label{conj z osc}
 All zero-energy orbits avoiding the $z$-axis have $z(t)$ oscillatory. 
\end{conjecture}

By oscillatory we mean that $z(t)$ has infinitely many zeros without being identically zero on any interval, which of course implies the existence of three consecutive zeros needed for our theorem.  We tried and failed to prove this conjecture using Sturm comparison methods.  It is necessary that we avoid the $z$-axis in this conjecture; we have counterexamples otherwise.  But the $z$-axis is special in the Heisenberg group -- it is the cut and conjugate locus for the origin -- and needs to be treated with care.  However, numerical evidence suggests that even those zero-energy orbits which do intersect the $z$-axis are self-similar, although perhaps in a different sense than those which do not intersect the $z$-axis, with $z(t)$ not necessarily oscillatory.  We provide details at the end of this section.  The upshot is the following.

\begin{conjecture}\label{all self sim}
All zero-energy orbits are self-similar.
\end{conjecture}


The self-similarity in Theorem \ref{self sim thm} is most clearly understood in new coordinates.  Let
\begin{align*}
s &= \frac{1}{4}\log ((x^2+y^2)^2+16z^2) \\
\theta &= \arg(x,y) \\ 
u &= \arg(x^2+y^2, 4z).
\end{align*}
Then the corresponding momenta $p_s$ and $p_{\theta}$ correspond precisely to the dilational and angular momenta; that is, $p_s=J$.  This allows us to compute the dilational and rotational factors $\lambda$ and  $\varphi$ for any given orbit satisfying our hypotheses as 
$$
\lambda = \exp (s(t_2) - s(t_0)) 
\qquad \text{and} \qquad
\varphi = \theta(t_2) - \theta(t_0).
$$
The rotational part is standard, but finding the correct coordinates for the dilational part requires solving a system of non-linear PDEs.
Armed with the dilational factor $\lambda$ we can now explicitly determine the time reparametrization $\tau$ as follows (if $\lambda=1$ then no reparametrization is necessary):
\begin{align*}
\xi(t) &=\frac{1}{2}\log_{\lambda}\left(1-(t-t_0)\frac{1-\lambda^2}{t_2-t_0}  \right)  \\
\tau(t) &= t_0+(t_2-t_0)\frac{1-\lambda^{2\xi(t) + 2\lfloor \xi(t) \rfloor}}{1-\lambda^2}.
\end{align*}

Suppose we are given a segment of a curve $c$ on a fundamental time domain $[t_0, t_2]$, and we want to know $c(t)$ at some future time $t$.
Then $c$ can be extended from $[t_0, t_2]$ to next fundamental domain $[t_2, t_2 + \lambda^2(t_2-t_0)]$ by appropriately rotating, dilating, and reparametrization time.  
We can then iterate to move forward to any future time (until possible collision), and can easily modify to move backward to past time.  
The function $\xi$ plays an important discretizing role:
its floor gives the number of iterations needed.
That is, $\lfloor \xi(t) \rfloor$ is the number of rotated and dilated copies of $c$ which need to be pasted together in order to reach $c(t)$.  
Note that consecutive ``periods" scale by $\lambda^2$, giving a geometric series, so we have the following corollary.

\begin{corollary}\label{strat cor}
Orbits satisfying the hypotheses of Theorem \ref{self sim thm} stratify according to Table \ref{table 2}. All collisions happen in finite time with 
$$t_{\text{col}}=t_0+\frac{t_2-t_0}{1-\lambda^2}.$$
\end{corollary}

\begin{table}
\caption{Summary of behavior}
\def\arraystretch{1.2}
  \begin{tabular}{ |  c | c | c | }
    \hline
    Dilational momentum & Dilation factor & Behavior of orbit \\ \hline 
    $J<0$ & $\lambda<1$ & Future collision, unbounded past \\ \hline 
    $J>0$ & $\lambda>1$ & Past collision, unbounded future \\ \hline
    $J=0$ & $\lambda=1$ & Periodic or quasi-periodic motion \\ \hline 
  \end{tabular}
  \label{table 2}
\end{table}

\begin{figure}
\begin{tabular}{ccc}
\includegraphics[height=.352\textheight]{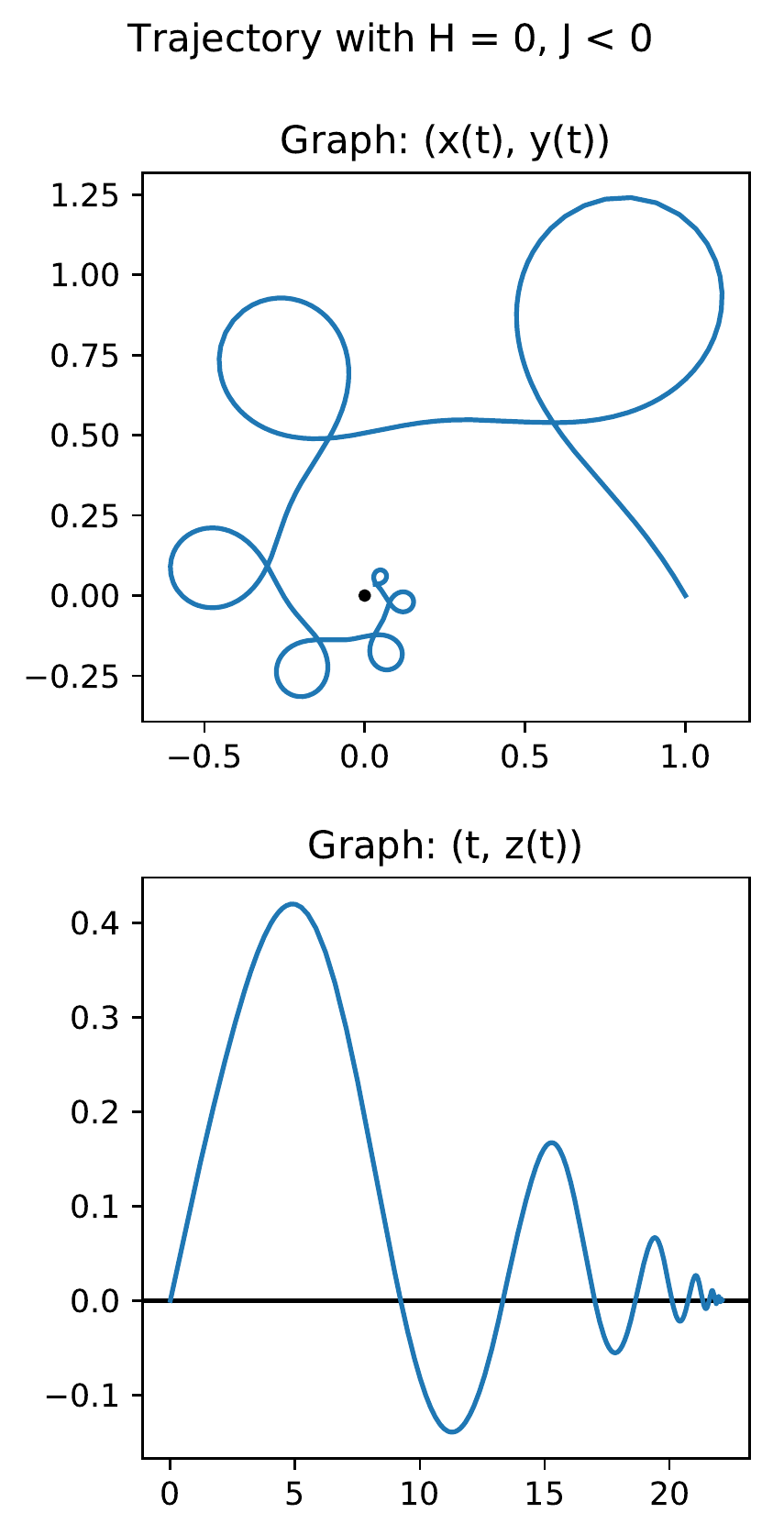}
\includegraphics[height=.352\textheight]{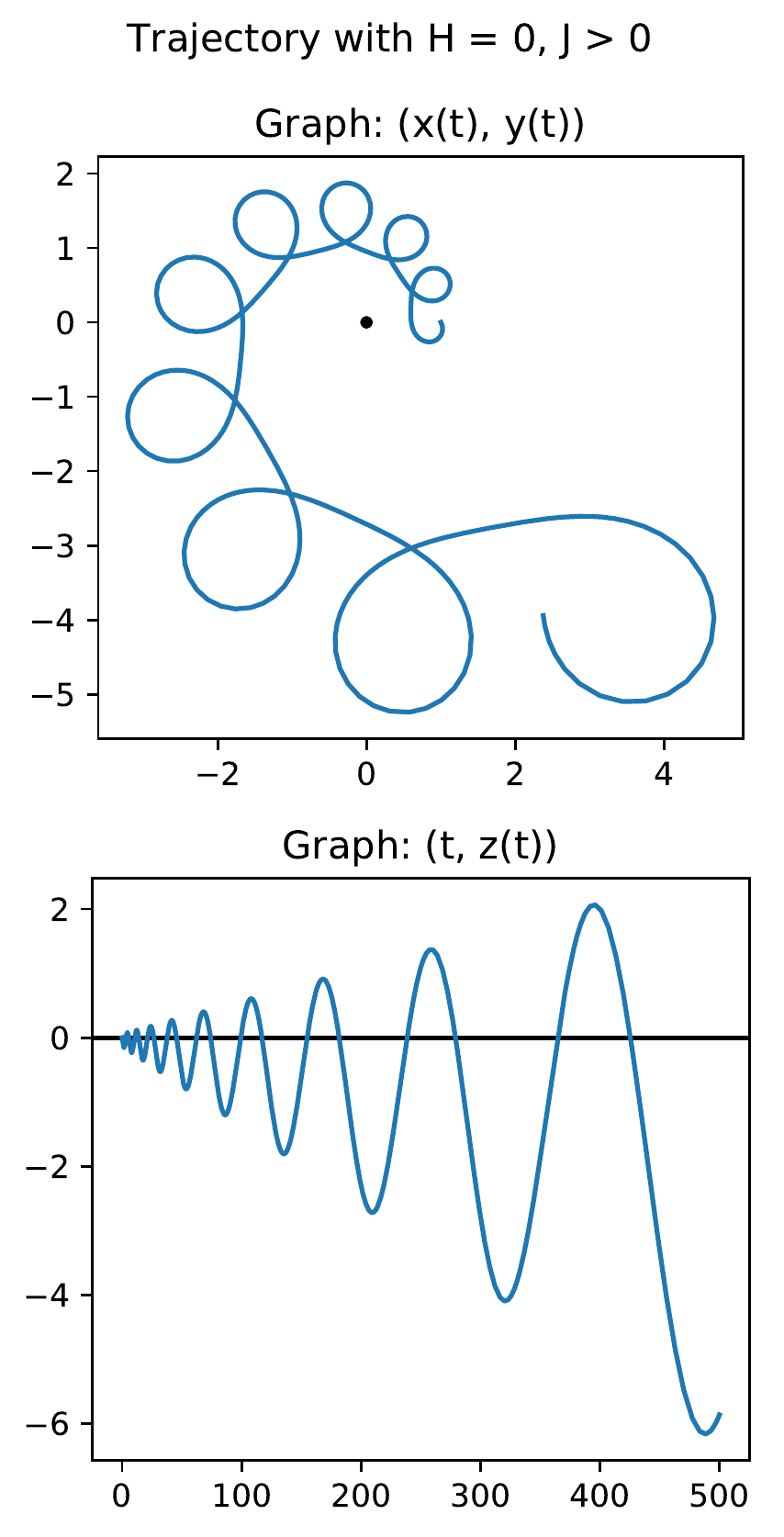} 
\includegraphics[height=.352\textheight]{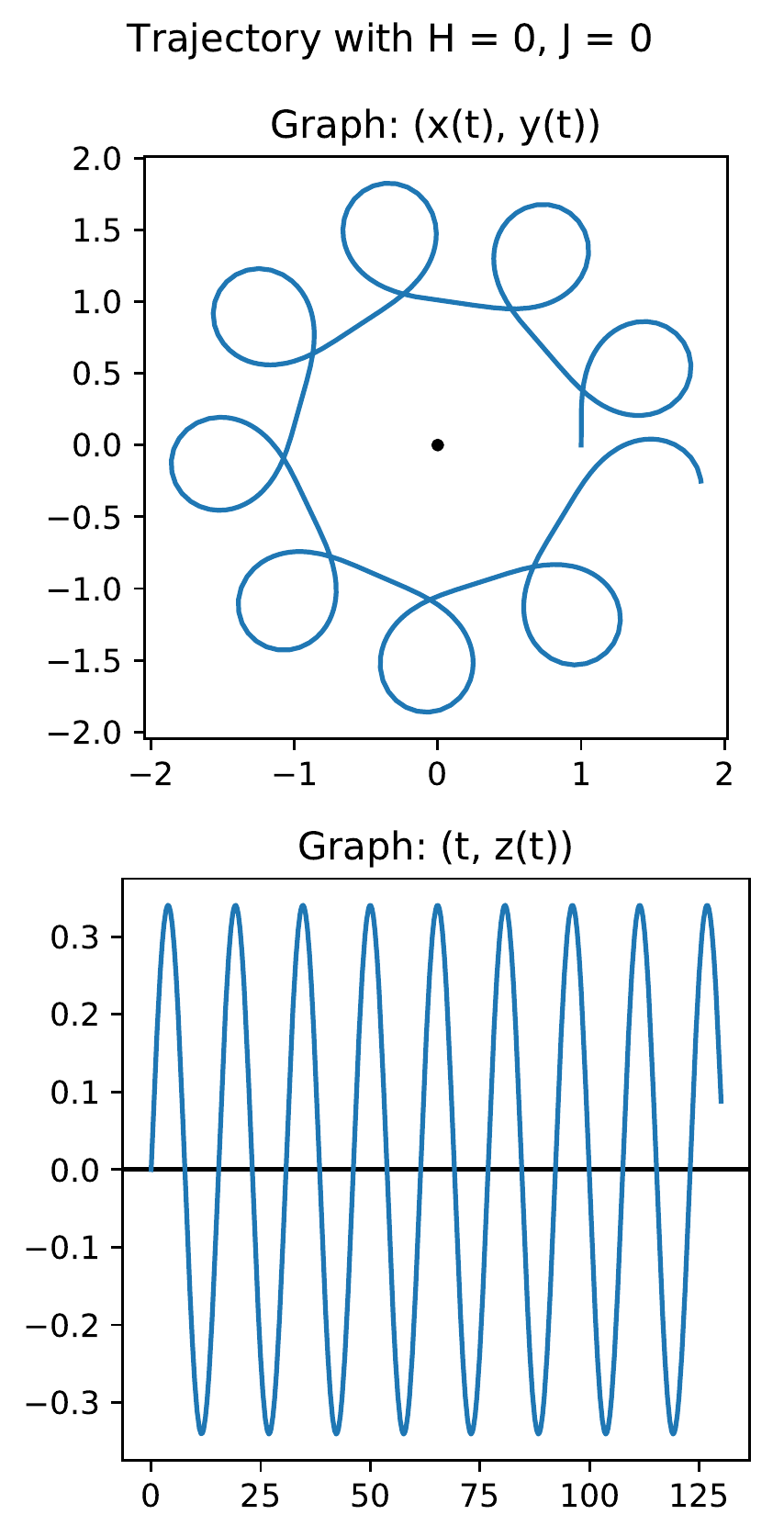}
\end{tabular}
\caption{Sample orbits with negative, positive, zero dilational momentum.} \label{strat fig}
\end{figure}

\begin{figure}
\centering 
\includegraphics[width=.9\textwidth]{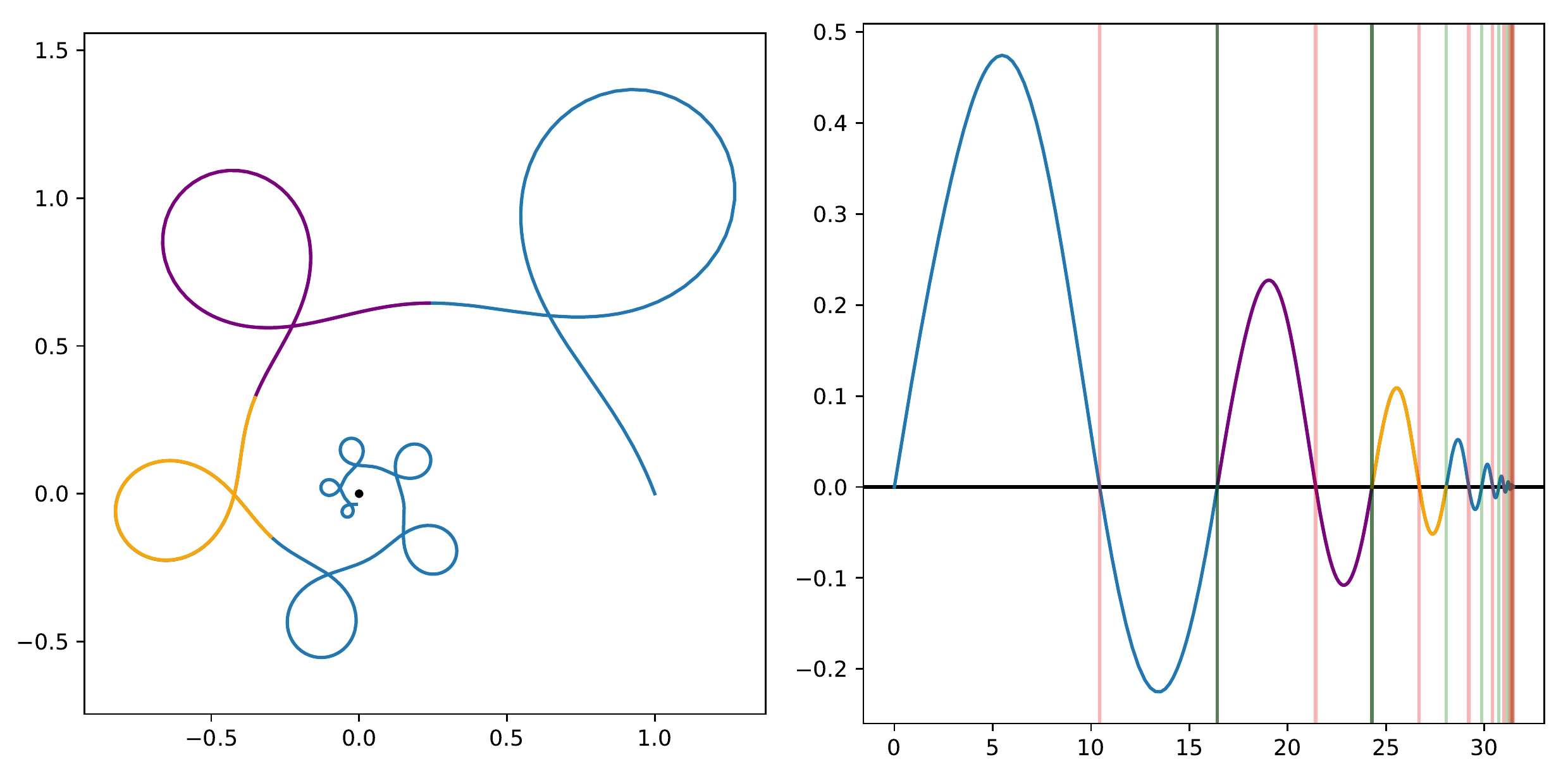}
\caption{A self-similar orbit.}\label{self sim fig}
\end{figure}

Figure \ref{strat fig} provides sample orbits from the stratification in 
Table \ref{table 2}.

Most of the self-similarity discussed in this section is apparent in Figure \ref{self sim fig}, which displays an orbit with negative dilational momentum.  On the left we see the projection of the orbit to the $xy$-plane, on the right we see $z$ as a function of time.  In both images the blue curve represents the orbit itself (numerically integrated), the purple segment represents the part over a chosen fundamental domain, and the orange segment is a rotated and dilated version of the purple.  The idea is that we are given the purple segment, from which we can recover the orange segment, which in turn matches the blue curve.  This can be iterated forward or backward in time.
The vertical lines on the right show zeros of $z(t)$; the given fundamental domain $[t_0, t_2]$ is shown with bold lines.  We can see how these ``periods" shrink geometrically and a collision occurs in finite future time.  
The codebase used to create Figure 
 \ref{self sim fig} will also be available at \cite{Github}.

To conclude this section we revisit Conjecture \ref{all self sim} and the discussion preceding it.   Recall that the $z$-axis is the conjugate and cut locus for the origin in the Heisenberg group. Zero-energy orbits meeting the $z$-axis need not have $z(t)$ oscillatory, but we believe they are all still self-similar.  However, numerical investigations suggest that such orbits bifurcate into two families with qualitatively different behavior, and this bifurcation occurs when dilational momentum $|J|=\frac{1}{2\sqrt{\pi}}$.

\begin{figure}[h]
\centering 
\includegraphics[width=.6\textwidth]{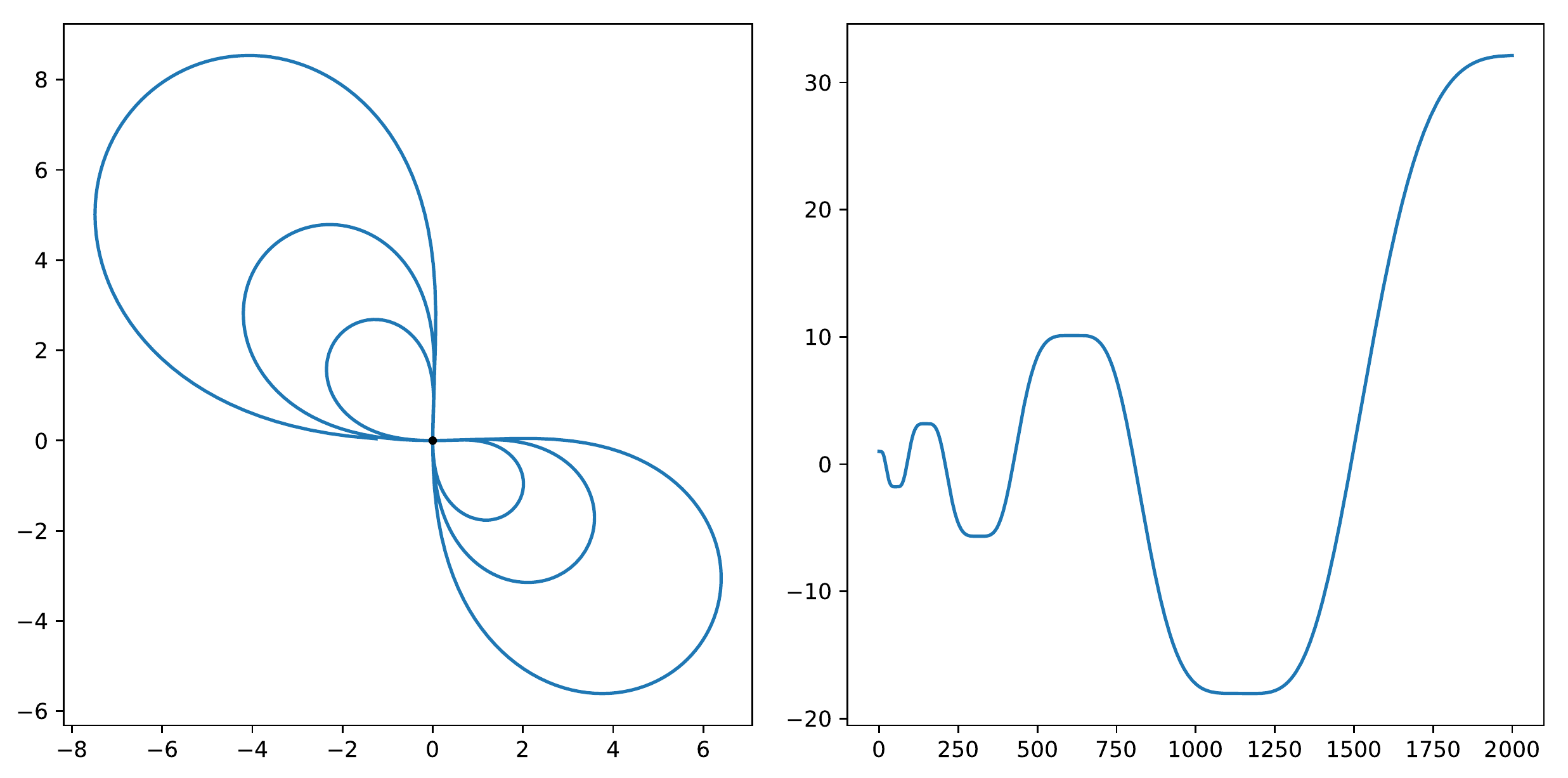}
\caption{Starting on $z$-axis with small dilational momentum.} \label{J small}
\end{figure}
\begin{figure}[h]
\centering 
\includegraphics[width=.6\textwidth]{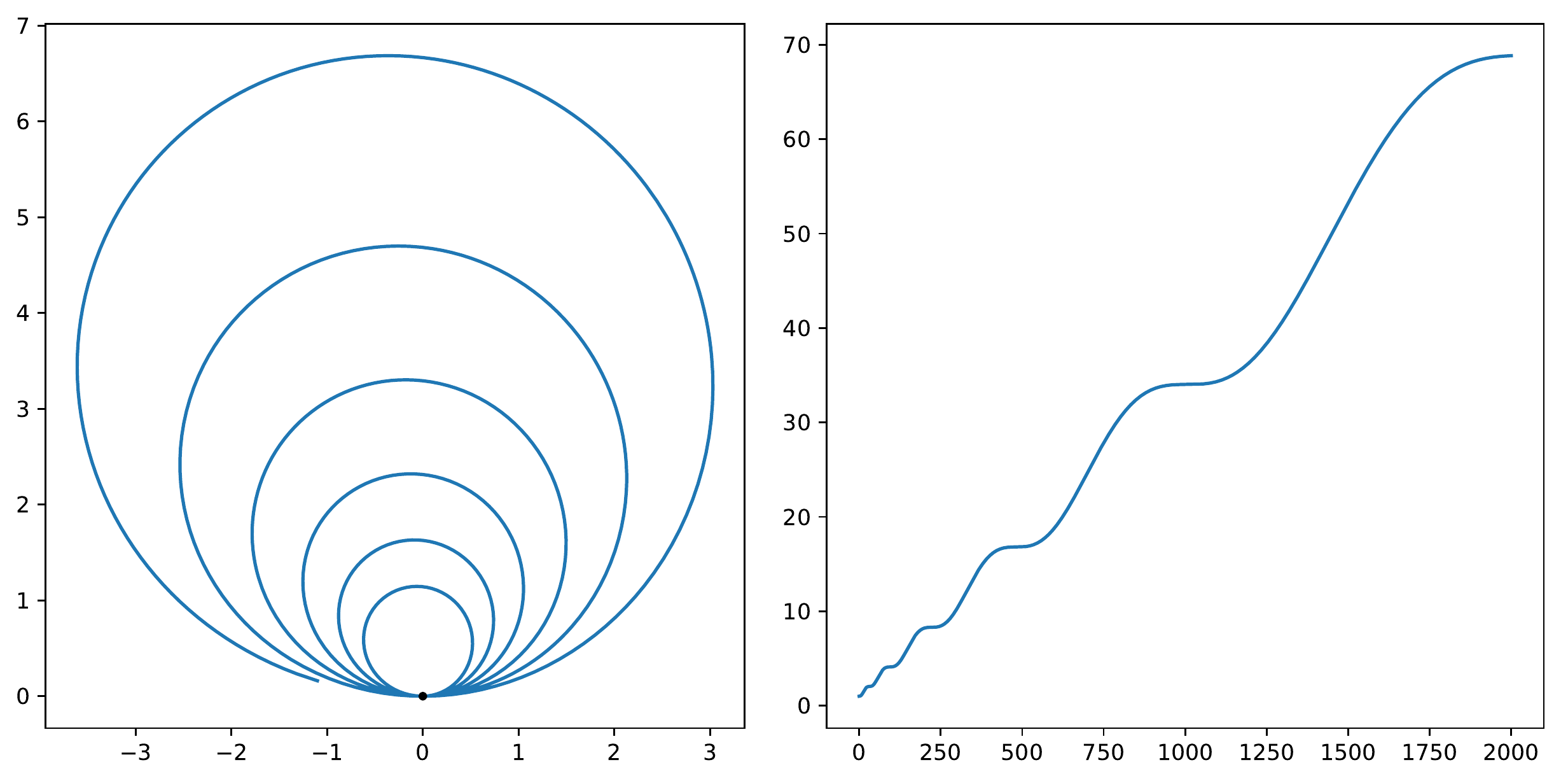}
\caption{Starting on $z$-axis with large dilational momentum} \label{J big}
\end{figure}

\begin{conjecture}  \label{conj z axis}
Consider a zero-energy orbit meeting the $z$-axis.
If $|J|<\frac{1}{2\sqrt{\pi}}$ then we have
self-similarity like a generic $H=0$ orbit. In particular, $z(t)$ oscillates and the orbit looks like a dilating figure eight.
If $|J|>\frac{1}{2\sqrt{\pi}}$ then we have
 self-similarity like Heisenberg geodesic.  In particular, $z(t)$ is monotonic and the orbit looks like a dilating helix. 
 See Figures \ref{J small} and \ref{J big}.
\end{conjecture}

\section{Open Questions} \label{open}
Here we simply list some of the many questions we have about this system and its generalizations.  Some of these we have investigated without success, others we have barely considered.

\begin{enumerate}
\item Is the Kepler-Heisenberg system integrable?  Are there any other non-trivial symmetries?  

\item How can we explicitly integrate the zero-energy system?  

\item How do orbits with nonzero energy behave?  Nearly all our efforts have focused on the $H=0$ case, where the dynamics are integrable and where the periodic orbits live.

\item Are Conjectures \ref{conj z osc},  \ref{all self sim}, \ref{conj z axis} true?

\item Can we regularize collisions with the sun?

\item Can we quantize this system in a meaningful way?

\item 
The richness of this system gives us hope that other ``sub-Riemannian mechanics" problems may be worth studying.  In each case we would take the sub-Riemannian Hamiltonian as our kinetic energy and choose some appropriate potential energy.  But this gives freedom to choose both the geometry and the potential:
\begin{enumerate}
 \item  Can we study the Kepler problem on higher dimensional Heisenberg groups?  Other Carnot groups?  Other sub-Riemannian spaces?  Is there a non-Euclidean setting where all three Kepler laws hold?
 \item Can we study the two body problem (not same as Kepler!) on the Heisenberg group?  Three bodies?   Other classical mechanics problems?  Other potentials?
 \item Combine the previous two items.  The $n$-body problem in jet space?
\end{enumerate}
\end{enumerate}


\begin{thebibliography}{99}

\bibitem{Diacu}  F. Diacu, E. Perez-Chavela, and M. Santoprete, The $n$-body problem in spaces of constant curvature, arXiv:0807.1747 (2008).

\bibitem{Github} V. Dods, ``Kepler-Heisenberg Problem Computational Tool Suite."  Available at \url{https://github.com/vdods/heisenberg}, 2018. 

\bibitem{DS} V. Dods and C. Shanbrom, Numerical methods and closed orbits in the Kepler-Heisenberg problem, \emph{Exp. Math.}, \textbf{28} (2019), 420--427.

\bibitem{DS2} V. Dods and C. Shanbrom, Self-similarity in the Kepler-Heisenberg problem, to appear in \textit{J. Nonlinear Sci.},
 arXiv:1912.12375.

\bibitem{Fiorani} E. Fiorani, G. Giachetta, and G. Sardanashvily, The Liouville-Arnold-Nekhoroshev theorem for non-compact invariant manifolds, \emph{J. Phys. A}, \textbf{36} (2003), L101-–L107.

\bibitem{Folland} G. Folland, A fundamental solution to a subelliptic operator,
\emph{Bull. Amer. Math. Soc.}, \textbf{79} (1973), 373--376.

\bibitem{Gromov} M. Gromov, \emph{Metric structures for Riemannian and Non-Riemannian Spaces}, Birkhauser (2007).

\bibitem{Katok} A. Katok and B. Hasselblatt, \emph{Introduction to the Modern Theory of Dynamical Systems}, Cambridge University Press (1995).

\bibitem{Tour} R. Montgomery, \emph{A Tour of Subriemannian Geometries, Their Geodesics and Applications}, AMS Mathematical Surveys and Monographs, \textbf{91} (2002).

\bibitem{MS} R. Montgomery and C. Shanbrom, Keplerian motion on the Heisenberg group and elsewhere, \emph{Geometry, mechanics, and dynamics: The Legacy of Jerry Marsden}, Fields Inst. Comm., \textbf{73} (2015), 319--342.

\bibitem{CS} C. Shanbrom, Periodic orbits in the Kepler-Heisenberg problem, \emph{J. Geom. Mech.}, \textbf{6} (2014), 261--278.

\bibitem{Tabachnikov} S. Tabachnikov, Geometry and Billiards, AMS Student Mathematical Library, \textbf{30} (2005).

\end{thebibliography}
\end{document}